\documentclass[12pt]{article}
\usepackage{amsfonts}
\usepackage{amsmath,amssymb} 
\newtheorem{thm}{Theorem}[section]

\newtheorem{cor}[thm]{Corollary}
\newtheorem{lem}[thm]{Lemma}
\newtheorem{propn}[thm]{Proposition}
\newtheorem{rem}[thm]{Remark}

\def\R{{\mathbb R}}
\def\E{{{\mathbb E}\,}}
\def\N{{\mathbb N}}
\def\P{{\mathbb P}}
\def\Z{{\mathbb Z}}

\def\sA{{\cal A}}
\def\sS{{\cal S}}

\def\ce{{\cal E}}
\def\sE{{\cal E}}
\def\cf{{\cal F}}
\def\sF{{\cal F}}

\def\proof{{\noindent{\sc Proof.}~}}

\def\qed{{\hfill $\square$ \bigskip}}

\newcommand{\dis}{\displaystyle}
\newcommand{\form}{{\cal E}}
\newcommand{\dom}{{\cal F}}
\newcommand{\rd}{{\mathbb R}^d}
\newcommand{\Sn}{{\mathcal S}_n}

\def\vareps {\varepsilon}
\def\eps{\varepsilon}

\def\fal{{~~~\mbox{for all }~}}
\def\tfrac#1#2{{\textstyle {#1\over #2}}}
\def\ol{\overline}

\def\lam{{\lambda}}

\def\norm#1{{\Vert #1 \Vert}}
\def\angel#1{{\langle#1\rangle}}

\makeatletter
\@addtoreset{equation}{section}

\makeatother
\begin{document}
\title{Convergence of symmetric Markov chains on $\Z^d$}
\author{Richard F. Bass, 
\thanks{Research partially supported by NSF grant DMS-0601783.}
 ~~Takashi Kumagai
\thanks{Research partially supported by the Grant-in-Aid for Scientific 
Research (B) 17340036 (Japan).} 
~ and ~
Toshihiro Uemura
\thanks{Research partially supported by the Grant-in-Aid for Scientific 
Research (C) 20540130 (Japan).}
}

\maketitle
\begin{abstract}
For each $n$ let  $Y^n_t$ be a continuous time symmetric Markov chain with state
space $n^{-1} \Z^d$. A condition in terms of
the conductances is given for the convergence
of the $Y^n_t$ to a symmetric Markov process $Y_t$ on $\R^d$. 
We have weak convergence of $\{ Y^n_t: t\leq t_0\}$ for every $t_0$
and every starting point.
The limit process $Y$ has a continuous
part and may also have jumps.
\end{abstract}

\section{Introduction}\label{S:intro}

For each $n$, let $Y^n_t$ be a continuous time symmetric Markov chain with
state space $\sS_n=n^{-1}\Z^d$ and
conductances $C^n(x,y)$. 
This means that $Y^n$ stays at a state $x$ for an exponential length of
time with parameter 
$\sum_{z\ne x} C^n(x,z)$ and then jumps to the next state $y$ with probability 
$C^n(x,y)/\sum_{z\ne x} C^n(x,z)$. It is natural to expect that one can give 
conditions on the conductances such that for each starting point and each 
$t_0$, the processes $\{Y^n_t; t\leq t_0\}$ converge weakly to a limiting 
process and that the limiting process be a symmetric Markov process. 
The purpose of this paper is to give such a theorem.

The earliest convergence theorem of this type is that of \cite{DFGW}
in the context of a central limit theorem for random walks in random
environment. A more general result is implicit in \cite{SZ}. In
\cite{BaKu08} the first two authors of the current paper 
extended the theorem in \cite{SZ}
in two ways: chains with unbounded range were allowed and the rather
stringent 
continuity conditions in \cite{SZ} were weakened. A chain with
unbounded range is one where there is no bound on the size of the jumps.
In all of these papers the limit process is a symmetric diffusion on
$\R^d$.

The paper \cite{HK07} considered conductances that were comparable to
the distribution  of a stable law and the limit process is what
is known as a stable-like process. Here the limit process has paths
that have no continuous part. A theorem for convergence of 
pure jump symmetric processes on $\R^d$ can be found in \cite{BKK}; as noted
there the methods can be readily modified to give a result on the 
convergence of symmetric Markov chains whose limiting process has
a more general jump structure than stable-like. Finally, we should
mention the well-known results of  \cite[Chap. 11]{SV} on 
non-symmetric Markov chains.

The current paper is devoted to proving a fairly general convergence
theorem for symmetric Markov chains. We point out three significant
differences from earlier work.
\begin{description}
\item{$\bullet$} Our Markov chains can have unbounded range and
the limit process is associated with a Dirichlet form with
both local and non-local components. This means the limit process
has a continuous part and may also have a discontinuous part.
\item{$\bullet$} We dispense with any continuity conditions on
the conductances. Instead only convergence locally in $L^1$ is needed.
\item{$\bullet$} The proofs are considerably simpler than previous work.
\end{description}

Let us give a heuristic description of our results, with the main
theorem stated precisely  in Section \ref{S:wc} as Theorem \ref{clt-main}. First of all,
the limiting symmetric Markov process is associated to the Dirichlet
form
$$\sE(f,f)=\int_{\R^d} \nabla f\cdot a\nabla f\, dx
+\int_{\R^d} \int_{\R^d}  (f(x)-f(y)^2 j(x,y)\, dx\, dy.$$
Here $a_{ij}(x)$ is a symmetric uniformly positive definite and bounded matrix
function. The first term on the right hand side represents the continuous
part of the limit process; if the second term on the right hand
side were not present, one would have a symmetric diffusion, and
the Dirichlet form would be the one arising from elliptic operators on $\R^d$
in divergence form. The double integral on the right hand side
represents the jump part, and very roughly says that the process
jumps from $x$ to $y$ with jump intensity $j(x,y)$.

We write our conductances as $C^n=C^n_C+C^n_J$, where $C^n_C$ and $C^n_J$
are the local (continuous)  and non-local (jump) parts, resp. Let us
discuss the local part first.
If one wants to
understand the behavior of the limiting process at a point $x$,
say, to look at $a(x)$,
a bit of thought leads to the realization that jumps by the Markov
chains that jump over but do not land on $x$ contribute. Thus, in
one dimension, one looks at a quantity $a^n(x)$  involving sums
of terms involving $C^n_C(y,z)$ with
$y\le x\le z$.
In higher dimensions one uses  a similar idea: one looks at the contribution
of $C^n_C(y,z)$ where $x$ lies on the shortest path from $y$ to $z$; a path
here means that at each step the path goes from a point to one of its
nearest neighbors. 
There is no single shortest path in general, so we form $a^n_{ij}(x)$
in terms of an average of expressions  involving $C^n_C(y,z)$, the
average being over all shortest paths from $y$ to $z$ that pass through
$x$. There are some very mild regularity conditions on $C^n$, but the
main hypothesis is that the $a^n_{ij}(x)$ are uniformly bounded and converge
to $a_{ij}(x)$ locally in $L^1$.

The conditions on the jump part are even weaker. We form a measure
$j^n(x,y)\, dx\, dy$ in terms of the $C^n_J$. We then require
that for each $N$, the measure $j^n(x,y)\, dx\, dy$ restricted to
$B_N=(B(0,N)\times B(0,N)) \setminus (B(0,N^{-1})\times B(0, N^{-1}))$
converges weakly to the measure $j(x,y)\, dx\, dy$ restricted to $B_N$,
where $B(0,r)$ is the ball of radius $r$ centered at 0.

After giving some definitions and setting up the framework in 
Section \ref{S:frame}, we obtain upper and lower bounds and
regularity results for the heat kernels for $Y^n$ in Sections
\ref{S:HKE} and \ref{S:lb}. The formulation of the main theorem
is given in Section \ref{S:wc} and the proof is given in Section 
\ref{S:proof}.

\section{Framework}\label{S:frame}

For $n\in \N$, let $\sS_n=n^{-1}\Z^d$.  Let $|\cdot|$
be the Euclidean norm and $B_n(x,r):=\{y\in \sS_n: |x-y|< r\}$. 

For $n\in \N$, let $C^n(\cdot,\cdot)$ be a symmetric function defined on
$(\sS_n\times\sS_n)\setminus\Delta$ into $\R_+$, where $\Delta=\{(x,x): x\in \sS_n\}$.  
Here symmetric means $C^n(x,y)=C^n(y,x)$ for all $x\ne y$.  
We call $C^n(x,y)$ the {\it conductance} between $x$ and $y$. 
Throughout the paper, we assume the following;

\medskip
\noindent (A1) {\it There exist $c_1,c_2>0$ independent of $n$ such that 
\[c_1\le \nu_x^n:=\sum_{y\in \sS_n}C^n(x,y)\le c_2\fal x\in \sS_n.\]
}

\medskip
\noindent (A2){\it There exist $M_0\ge 1,
\delta>0$ independent of $n$ such that the following holds: 
for any $x,y\in \sS_n$ with $|x-y|=n^{-1}$, there exist $N\ge 2$ and 
$x_1,\cdots, x_N\in B_n(x,n^{-1}M_0)$ such that $x_1=x$, $x_N=y$ 
and $C^n(x_i,x_{i+1})\ge \delta$ for $i=1,\cdots, N-1$. 
}

\medskip
\noindent (A3) {\it There exists a function $\varphi : {\mathbb R}_+ 
\longrightarrow {\mathbb R}_+$ so that for any $n\in {\mathbb N}$,
$$
C^n(x,y) \le  n^{-(d+2)} \varphi\bigl(|x-y|\bigr), \quad x,y\in \sS_n 
\quad {\rm and} \quad 
\dis{\int_0^{\infty} \bigl(1\wedge t^2\bigr) \, t^{d-1}
\varphi(t)dt<\infty}.
$$
}

\medskip
Note that from the assumption (A3), we see for any $x\in \sS_n$,
$$
\begin{array}{lrl}
\dis{n^2 \sum_{y\in \sS_n} \bigl(1\wedge |y|^2\bigr) C^n(x,x+y)} 
&\le& \dis{ n^2 \sum_{y\in \sS_n} \bigl(1\wedge |y|^2\bigr) 
n^{-(d+2)} \varphi(|y|) =\sum_{y\in \sS_n} \bigl(1\wedge |y|^2\bigr)
\varphi(|y|) n^{-d} } \\

& & \vspace{-6pt} \\
&\le& c_d\dis{\int_{\rd} \bigl(1\wedge |y|^2\bigr) \varphi(|y|) dy 
=c_d' \int_0^{\infty} (1\wedge t^2)\,t^{d-1} \varphi(t)dt<\infty.}
\end{array} 
$$
Thus we have 
\begin{equation}
\dis{
M:=\sup_n \sup_{x\in \sS_n} \Bigl( 
n^2 \sum_{y\in \sS_n} \bigl(1\wedge |y|^2\bigr) C^n(x,x+y)} \Bigr) 
<\infty.
\label{cond-bound}
\end{equation}

An example of $C^n(x,y)$ that satisfies (A1), (A2) and (A3) is the following: 
\begin{align*}
\frac {c_11_{\{1\geq |x-y|\ge n^{-1}\}}}{n^{d+2}|x-y|^{d+\alpha}}
&+c_21_{\{|x-y|= n^{-1}\}}\le C^n(x,y)\\
&\le 
\frac {c_31_{\{1\geq |x-y|\ge n^{-1}\}}}{n^{d+2}|x-y|^{d+\beta}}
+c_41_{\{|x-y|= n^{-1}\}}
+\frac{c_5 1_{\{|x-y|>1\}}} {n^{d+2}|x-y|^{d+\alpha}},
\end{align*}
where $0<\alpha\le \beta<2$. 

Let $\mu^n_x\equiv n^{-d}$ for all $x\in \sS_n$ and for each $A\subset \sS_n$, define 
$\mu^n (A)=\sum_{y\in A}\mu^n_y$ and $\nu^n (A)=
\sum_{y\in A}\nu_y^n$. 
Note that $L^2(\sS_n, \mu^n)=L^2(\sS_n, \nu^n)$
by (A1). Now, for each $f\in L^2(\sS_n,\mu^n)$, define 
\begin{eqnarray}
\sE^n (f,f)&=&\frac{n^{2-d}}2 \sum_{x,y\in \sS_n} (f(y)-f(x))^2 C^n(x,y),\label{DFe_n}\\
\cf^n&=&\{f\in L^2(\sS_n,\mu^n): \ce^n(f,f)<\infty\}.\label{DFf_n}
\end{eqnarray}
For $p\ge 1$, define $\|f\|_{p,n}^p=\sum_{y\in \sS_n} \bigl|f(x)\bigr|^p 
\mu^n_x$. The following lemma is standard. 
\begin{lem}\label{dfl^2} For each $n$, $\dom^n=L^2(\Sn,\mu^n)$, and for 
$f\in L^2(\Sn,\mu^n)$, we have 
$$
\form^n(f,f) \le 2 n^2 M ||f||_{2,n}^2,
$$
where $M$ is the constant appearing in (\ref{cond-bound}). 
\end{lem}
\proof  Let $f\in L^2(\Sn, \mu^n)$. Since $|x-y|\ge 1/n$ 
for any $x,y\in\Sn$ with $x\not=y$, we have  
$$
\begin{array}{lcl}
\dis{ \frac{n^{2-d}}2 \sum_{x,y\in \Sn \atop x\not= y} \bigl(f(x)-f(y)\bigr)^2 
C^n(x,y) } &\le& \dis{ n^{2-d} \sum_{x,y\in \Sn \atop |x-y|\ge 1/n} 
\bigl(f(x)^2+f(y)^2\bigr)  C^n(x,y)  } \\
& & \vspace{-6pt} \\
&\le& \dis{ 2 n^{2-d} \sum_{x\in \Sn} f(x)^2 
\sum_{y\in \Sn \atop |x-y|\ge 1/n} C^n(x,y) } \\
& & \vspace{-6pt} \\
&\le& \dis{ 2 n^{4-d} \sum_{x\in \Sn} f(x)^2 
\sum_{y\in \Sn \atop |x-y|\ge 1/n} \bigl( 1\wedge |x-y|^2\bigr) C^n(x,y) } \\
& & \vspace{-6pt} \\
&\le& \dis{ 2 n^{2-d} M \sum_{x\in \Sn} f(x)^2 =2 n^2 M 
||f||_{2,n}^2. } \\
\end{array}
$$\qed 

Using Lemma \ref{dfl^2}, it is easy to check that $(\sE^n,\sF^n)$ 
is a regular Dirichlet form on  $L^2(\sS_n,\mu^n)$. Further, $\sF^n=L^2(\sS_n,\mu^n)$ 
is equal to the 
closure of the space of compactly supported functions on 
$\sS_n$ with respect to $(\sE^n(\cdot, \cdot)+||\cdot||^2_{2,n})^{1/2}$. Let $Y^{(n)}_t$ be the
corresponding continuous time Markov chains on $\sS_n$ and let 
$p^n(t,x,y)$ be the transition density for $Y^{(n)}_t$ with 
respect to $\mu^n$.
The infinitesimal generator of $Y^{(n)}_t$ can be written as 
\[ \sA^n f(x)=\sum_{y\in \sS_n} (f(y)-f(x)) 
C^n(x,y)n^2
=\sum_{y\in \sS_n}(f(y)-f(x))\frac{C^n(x,y)n^{2-d}}{\mu^n_x},\]
for each $f\in L^2(\sS_n,\mu^n)$. 
\begin{rem}\label{thm:conserem}
{\rm 
Note that under (A1), $\{Y^{(n)}_t\}$ is conservative. Indeed, define  
a symmetric Markov chain $\{X^{(n)}_m\}$ by
\[\P^x(X^{(n)}_1=y)=\frac {C^n(x,y)}{\nu^n_x}\fal 
x,y\in \sS_n.\]
Then the corresponding semigroup satisfies $P^{X,n}_11(x)=\sum_{y\in \sS_n}
\P^x(X^{(n)}_1=y)=1$ by (A1), so inductively we have $P^{X,n}_m1=1$ for all $m\in \N$,
so that $\{X^{(n)}_m\}$ is conservative. But $\{Y^{(n)}_t\}$ is a time changed process of 
$\{X^{(n)}_m\}$. To see this, 
let $\{U_i^{x,n}: i\in \N, x\in \sS_n\}$ be an independent sequence of exponential random 
variables, where the  parameter for $U_i^{x,n}$ is $\nu^n_x$, 
  that is independent of $X^{(n)}_m$,
and define $T^{(n)}_0=0, T^{(n)}_m=\sum_{k=1}^m U_k^{X^{(n)}_{k-1},n}$. 
Set $\widetilde Y^{(n)}_t=X_m$ if $T^{(n)}_m\le t<T^{(n)}_{m+1}$; then 
the laws of $\widetilde Y^{(n)}$ and $Y^{(n)}$ are the
same, and hence $\widetilde Y^{(n)}$ is a realization of the continuous time Markov
chain corresponding to (a time change of) $X^{(n)}_m$. 
Note that by (A1), the mean exponential holding time at each point 
for $\widetilde Y^{(n)}$ can be controlled uniformly from above and below by a positive 
constant, so we conclude $P^{n}_t1=1$ for all $t>0$, where $P^n_t$ is the semigroup 
corresponding to $\{Y^{(n)}_t\}$. 
}
\end{rem}

\section{Heat kernel estimates}\label{S:HKE}
\subsection{Nash inequality}
For $f\in L^2(\sS_n, \mu^n)$, let
\[\ce^n_{NN}(f,f)=\frac{n^{2-d}}2\sum_{{x,y\in \sS_n} \atop 
{|x-y|=n^{-1}}}(f(x)-f(y))^2,\]
which is the Dirichlet form for the simple symmetric random walk in $\sS_n$. 
By \cite[Proposition 3.1]{BaKu08}
there exists $c_1>0$ independent of $n$ such that for any $f\in L^2(\sS_n, 
\mu^n)$, 
\begin{equation}\label{eq:nash1}
\|f\|_{2,n}^{2(1+2/d)}\le c_1 \ce^n(f,f)\|f\|_{1,n}^{4/d}, 
\end{equation} 
and
\begin{equation}
p^n(t,x,y)\le c_1t^{-d/2}~~\fal x,y\in\sS_n, t>0.\label{eq:nash3}
\end{equation}

For $r\in (n^{-1},1]$, let $\sE^{n,r}$ be the Dirichlet form corresponding
to $\{Y^{(n),r}_t:=r^{-1}Y^{(n)}_{r^2t}, t\geq 0\}$.
By simple computations, we have 
\[\sE^{n,r} (f,f)=\frac{(nr)^{2-d}}2 \sum_{x,y\in \sS_{nr}} (f(y)-f(x))^2 C^n(rx,ry),\]
where $\sS_{nr}=\{x/r: x\in \sS_n\}=(nr)^{-1}\Z^d$. Define
\begin{equation}\label{eqn:003.8}
p^{n,r}(t, x, y):=r^{d}p^n(r^2t, rx, ry).
\end{equation}
Then $p^{n,r}(t,x,y)$ is the heat kernel for $\sE^{n,r}$. By (\ref{eq:nash3}), 
we have
\begin{equation}
p^{n,r}(t,x,y)\le c_1t^{-d/2}~~\fal x,y\in\sS_{nr}, t>0.\label{eq:nash4}
\end{equation}

For $\lambda\ge 1$, let $Y^{(n),r,\lambda}_t$ be a process on $\sS_{nr}$ 
with the  large jumps of $Y^{(n)}_t$ removed. More precisely, 
$Y^{(n),r,\lambda}_t$ is a process whose Dirichlet form is 
\[\ce^{n,r,\lambda} (f,f)=\frac 12 \sum_{{x,y\in \sS_{nr}}\atop
{|x-y|\le\lambda}}(f(x)-f(y))^2 (nr)^{2-d}C^n(rx,ry),\]
for each $f\in L^2(\sS_{nr},\mu^{nr})$. We denote the heat kernel for 
$Y^{(n),r,\lambda}_t$ by $p^{n,r,\lambda}(t,x,y)$, $x,y\in \sS_{nr}$.

\subsection{Exit time  probability estimates}

In this subsection, we will obtain some  exit time 
estimates. Note that similar estimates are obtained in \cite[Proposition 3.7]{Foo} and 
\cite{CK09}. 
\begin{propn}\label{thm:4.2}
For $A> 0$ and $0<B<1$, there exists  
$t_0 =t_0(A,B)\in (0,1)$ such that for every $n\in \N$, $r\in (0,1]$ and 
$x\in \sS_n$,
\begin{equation}\label{eq:tighty}
\P^x \left( \sup_{s\le r^2t_0}|Y^{(n)}_t-Y^{(n)}_0|>rA\right)
=\P^x \left( \sup_{s\le t_0}|Y^{(n),r}_t-Y^{(n),r}_0|>A\right)
\le B.\end{equation}
\end{propn}
\proof
Let $\lambda>0$. Since we have (\ref{eq:nash4}) and $p^{n,r,\lambda}(t,x,y)
\leq p^{n,r}(t,x,y)$, by Theorem (3.25) of \cite{CKS}, we have 
\begin{equation}\label{kern4}
p^{n,r,\lambda}(t,x,y)\leq c_1\;t^{-\frac{d}2}
\;\exp\left(-E(2t,x,y)\right)
\end{equation}
for all $t\le 1$ and $x,y\in \sS_{nr}$, where  
\begin{eqnarray*}
E(t,x,y)&=&\sup\{|\psi(y)-\psi(x)|-t\; \Lambda(\psi)^2 :
\Lambda(\psi)<\infty\},\\
\Lambda(\psi)^2&=& \|e^{-2\psi}\Gamma_{\lambda,r}[e^\psi]\|_\infty \vee
\|e^{2\psi}\Gamma_{\lambda,r}[e^{-\psi}]\|_\infty,
\end{eqnarray*}
and $\Gamma_{\lambda,r}$ is defined by 
\begin{equation}  \label{densi}
\Gamma_{\lambda,r}[v](\xi)=
\sum_{{\eta,\xi\in\sS_{nr}}\atop
{|\xi-\eta|\le\lambda}}(v(\eta)-v(\xi))^2
C^n(r\eta,r\xi)(nr)^2,\qquad 
\xi\in\sS_{nr}.
\end{equation}
Now let $R=|x-y|$ and let $\psi(\xi)=s(|\xi-x|\wedge R)$.
Then, $|\psi(\eta)-\psi(\xi)|\le s|\eta-\xi|$,
so that 
\[(e^{\psi(\eta)-\psi(\xi)}-1)^2\le |\psi(\eta)-\psi(\xi)|^2
e^{2|\psi(\eta)-\psi(\xi)|}\le cs^2|\eta-\xi|^2e^{2|\psi(\eta)-\psi(\xi)|}\]
for $\eta,\xi\in \sS_{nr}$ where $|\eta-\xi|\le \lambda$. Hence 
\begin{eqnarray*}
e^{-2\psi(\xi)}\Gamma_{\lambda,r}[e^\psi](\xi)&=&
\sum_{{\eta\in\sS_{nr}}\atop
{|\xi-\eta|\le\lambda}}(e^{\psi(\eta)-\psi(\xi)}-1)^2
{C^n(r\eta,r\xi)(nr)^2}\\&\le&
c_1s^2e^{2s\lambda}\sum_{{\eta\in\sS_{nr}}\atop
{|\xi-\eta|\le\lambda}}|\eta-\xi|^2
{C^n(r\eta,r\xi)(nr)^2}\\&=&
c_1s^2e^{2s\lambda}\sum_{{\eta'\in\sS_n}\atop
{|\xi'-\eta'|\le\lambda r}}|\eta'-\xi'|^2
{C^n(\eta',\xi')n^2}
\\&\le &
c_1s^2e^{2s\lambda}\Big(\sum_{{\eta'\in\sS_n}\atop
{|\xi'-\eta'|\le 1}}|\eta'-\xi'|^2
{C^n(\eta',\xi')n^2}+((\lambda r)^2\vee 1)\sum_{{\eta'\in\sS_n}\atop
{|\xi'-\eta'|\ge 1}}
{C^n(\eta',\xi')n^2}\Big)\\
&\le & c_2(\lambda^2\vee 1)s^2e^{2s\lambda}\le c_3e^{3s\lambda}(1+1/\lambda^2)
\end{eqnarray*}
for all $\xi\in \sS_{nr}$ where (A3) and $r\le 1$ are 
used in the third inequality.
We have the same bound when $\psi$ is replaced by $-\psi$,
so $\Lambda(\psi)^2\leq c_3e^{3s\lambda}(1+1/\lambda^2)$. Now, let 
$\lambda=A/(6d)$, $t_0\le 1\wedge\lambda^4=1\wedge (A^4/(6d)^4)$ 
and $s=(3\lambda)^{-1}\log (1/t^{1/2})>0$. Then, for each $t\le t_0$ and $R\ge A$, 
\begin{eqnarray}
p^{n,r,\lambda}(t,x,y)&\leq &c_4t^{-\frac{d}2}
\;\exp\left(-sR+c_3te^{3s\lambda}(1+1/\lambda^2)\right)\nonumber\\&\le&
c_5\exp\Big((d-\frac {2Rd}{A})\log (\frac 1{t^{1/2}})\Big)
\le c_5\exp\Big(-\frac {Rd}{A}\log (\frac 1{t^{1/2}})\Big).\label{tightwd1}
\end{eqnarray}
Thus, 
\begin{eqnarray}
\sum_{B_{nr}(x,A)^c}p^{n,r,\lambda}(t,x,y)\mu_y^{nr}
&\le & c\int_A^\infty R^{d-1}
\exp\Big(-\frac{Rd}{A}\log (\frac 1{t^{1/2}})\Big)dR\nonumber\\
&= & cA^d\int_1^\infty R'^{d-1}
\exp\Big(-R'd\log (\frac 1{t^{1/2}})\Big)dR'<B/4\label{tightwd2}
\end{eqnarray}
for all $t\le t_0$ if we choose $t_0$ small, depending on $A$ and $B$.
Thus, applying \cite[Lemma 3.8]{BBCK}, we obtain
\begin{equation}\label{lamkil111}
\P^x\left( \sup_{s\le t_0}|Y^{(n),r,\lambda}_t-Y^{(n),r,\lambda}_0|
>A\right)\le B/2.
\end{equation}
We now use Meyer's argument to obtain the estimate for $Y^{(n),r}$. 
Note that for any $x\in \sS_{nr}$, 
\begin{eqnarray*}
{\cal J}(x)&:=&\sum_{{y\in\sS_{nr}}\atop
{|x-y|\ge\lambda}}C^n(rx,ry)(nr)^2
\le \sum_{{y\in\sS_{nr}}\atop{|x-y|\ge\lambda}}
\frac{(r^2|x-y|^2)\wedge 1}{\lambda^2r^2}C^n(rx,ry)(nr)^2\\
&=& \frac {1}{\lambda^2} \sum_{y'\in \sS_n}\Big(|x'-y'|^2\wedge 1 \Big)
C^n(x',y')n^2\le \frac {M}{\lambda^2}=\frac {(6d)^2M}{A^2},
\end{eqnarray*}
where (A3) is used in the last inequality. So, if we let
$U_1:=\inf\{t>0: \int_0^t{\cal J}(Y^{(n),r}_s)ds>S_1\}$,
where $S_1$ is the independent exponential distribution with mean $1$, 
we have 
\begin{equation}\label{lamkil133}
P(U_1\le t_0)\le 1-e^{-(6d)^2t_0/A^2}<B/2
\end{equation} 
by taking $t_0$ small.
Using Meyer's argument (see, for example, Section 4.1 in \cite{CK08}), we obtain
\begin{align*}
\P^x\Big( \sup_{s\le t_0}|Y^{(n),r}_t-Y^{(n),r}_0|>A\Big)
&=\P^x\Big( \sup_{s\le t_0}|Y^{(n),r}_t-Y^{(n),r}_0|>A, U_1>t_0\Big)\\
&\qquad \qquad +\P^x\Big( \sup_{s\le t_0}|Y^{(n),r}_t-Y^{(n),r}_0|>A, U_1\leq t_0\Big)\\
&\le  \P^x\left( \sup_{s\le t_0}|Y^{(n),r,\lambda}_t-Y^{(n),r,\lambda}_0|>A\right)
+\P^x\left(U_1\le t_0\right)\\
&\le B/2+B/2=B,
\end{align*}
where  (\ref{lamkil111}) and (\ref{lamkil133}) are used in the last inequality.
\qed
\begin{cor}\label{coe:4.2q1}
For $0<A', B'<1$, there exists 
$R_0 =R_0(A',B')>0$, such that for every $n\in \N$, $r\in (0,1]$ and 
$x\in \sS_n$,
\begin{equation}\label{eq:tighwty12}
\P^x \left( \sup_{s\le r^2A'}|Y^{(n)}_t-Y^{(n)}_0|>rR_0\right)
=\P^x \left( \sup_{s\le A'}|Y^{(n),r}_t-Y^{(n),r}_0|>R_0\right)
\le B'.\end{equation}
\end{cor}
\proof In the proof of Proposition \ref{thm:4.2}, 
take $A\ge 1$, $\lambda=A^{1/2}/(6d)$ (instead of $\lambda=A/(6d)$) and
$A\ge 1$. Then, since $A^{1/2}\le A\le R$, we have (\ref{tightwd1}) 
by changing $A$ to $A^{1/2}$. So as in (\ref{tightwd2}), there exists 
$R_0$ large such that for $t\le t_0=:A'$ and $A\ge R_0$, we have 
\begin{eqnarray*}
&&\sum_{B_{nr}(x,A)^c}p^{n,r,\lambda}(t,x,y)\mu_y^{nr}
\le  cA^{d/2}\int_{A^{1/2}}^\infty R'^{d-1}
\exp\Big(-R'd\log (\frac 1{t^{1/2}})\Big)dR'\\
&\le & cA^{d/2}\exp\Big(-\frac {A^{1/2}}2 d\log (\frac 1{t^{1/2}})\Big)
\int_{A^{1/2}}^\infty R'^{d-1}
\exp\Big(-\frac {R'}2 d\log (\frac 1{t^{1/2}})\Big)dR'<B/4.
\end{eqnarray*} 
Also, similarly to (\ref{lamkil133}), we have 
\[P(U_1\le t_0)\le 1-e^{-(6d)^2t_0/A}<B/2\]
for all $A\ge R_0$, by taking $R_0$ large. With these changes, we can 
obtain the result similarly to the proof of Proposition \ref{thm:4.2}.\qed

\section{Lower bounds and regularity for the heat kernel}\label{S:lb}
We now introduce the space-time process
$Z^{(n)}_s:=(U_s, Y^{(n)}_s)$, where $U_s=U_0+s$. The filtration generated
by $Z^{(n)}$ satisfying the usual conditions
will be denoted by $\{ \widetilde \cf_s; \, s\geq 0\}$.
The law of the space-time
process $s\mapsto Z^{(n)}_s$ starting from $(t, x)$ will be denoted 
by $\P^{(t, x)}$.
We say that  a non-negative Borel  measurable function
$q(t,x)$ on $[0, \infty)\times \sS_n$ is {\it parabolic}
in a relatively open subset $B$ of $[0, \infty)\times \sS_n$ if
for every relatively compact open subset $B_1$ of $B$,
$q(t, x)=\E^{(t,x)} \left[ q (Z^{(n)}_{\tau^n_{B_1}}) \right]$
for every $(t, x)\in B_1$,
where $\tau^n_{B_1}=\inf\{s> 0: \, Z^{(n)}_s\notin B_1\}$.

\medskip

We denote $T_0:=t_0(1/2, 1/2)<1$ the constant
in (\ref{eq:tighty}) corresponding to $A=B=1/2$.
For $t\ge 0$ and $r>0$, we define
\[  Q^n(t,x,r):=[t,t+T_0 r^{2}]\times B_n(x,r),
\]
where $B_n(x,r)=\{y\in \sS_n: |x-y|< r\}$. 

It is easy to see the following 
(see, for example, Lemma 4.5 in \cite{CK} for the proof).

\begin{lem}\label{4.6}
For each $t_0>0$ and $x_0\in \sS_n$, $q^{n}(t,x):=p^{n} (t_0-t, x,x_0)$
is parabolic on $[0,t_0)\times \sS_n$.
\end{lem}

For $A\subset \sS_n$ and a process $Z_t$ on $\sS_n$, let 
\[ 
\tau^n=\tau^n_A(Z):=\inf\{t\ge 0:
Z_t\notin A\}, \qquad T^n_A=T^n_A(Z):=\inf\{t\ge 0: Z_t\in A\}. \]

The next proposition provides a lower bound for the heat kernel
and is the key step for the proof of the H\"older
continuity of $p^n (t, x, y)$. 

\begin{propn}\label{szhol} 
There exist $c_1>0$ and $\theta\in (0,1)$ such that 
for each $n\in \N$, if $|x-x_0|,|y-x_0|\leq t^{1/2}$,
$x,y,x_0\in \sS_n$, $t\in (n^{-1},1]$  and $r\geq t^{1/2}/\theta$, then
\[ \P^x(Y^{(n)}_t=y, \tau_{B(x_0,r)}^n>t)\geq c_1t^{-d/2}n^{-d}. \]
\end{propn}

To prove this we first need some preliminary lemmas.
The proof of the following weighted Poincar\'e inequality can be found in 
\cite[Lemma 1.19]{SZ} and \cite[Lemma 4.3]{BaKu08}.

\begin{lem}\label{1.19}
Let
\[ 
g_n(x)= c_1 \prod_{i=1}^d e^{-|x_i|}\qquad x\in \sS_n,
\]
where $c_1$ is determined by the equation 
$\sum_{l\in \sS_n}g_n(x)\mu^n_x=n^d$. Then there exists $c_2>0$ such that 
\[
c_2\Big<(f-\langle f \rangle_{g_n})^2\Big>_{g_n}\le n^{2-d} 
\sum_{l\in \sS_n}g_n(l) \sum_{i=1}^d\Big(f(l+\frac{e^i}n)-f(l)\Big)^2, 
\qquad  f\in L^2(\sS_n), 
\]
where 
\[
\langle f \rangle_{g_n}=\sum_{l\in \sS_n}f(l)g_n(l)\mu_l^n
\]
and $e^i$ is the element of $\Z^d$ whose $j$-th component is $1$ if $j=i$ 
and $0$ otherwise. 
\end{lem}

We now give a key lemma. 
\begin{lem}\label{1.20}
There is  an $\eps>0$ such that 
\begin{equation}\label{1.16}
p^n(t,x,y)\ge \eps t^{-d/2},
\end{equation}
for all $n\in \N$, $(t,x,y)\in (n^{-1},1]\times \sS_n\times \sS_n$ with 
$|x-y|\le 2t^{1/2}$.
\end{lem}
\proof 
It is enough to prove the following:
there is an $\eps>0$ such that 
\begin{equation}\label{1.25}
(nr)^{-d}\sum_{l\in \sS_{nr}}\log \Big(p^{n,r}(\tfrac 12, k,l+m)\Big) 
g_{nr}(l)\ge \tfrac 12 \log\eps,
\end{equation} 
for any $n\in \N$, $r\in (n^{-1},1]$ and $k,m\in \sS_n$ with $|k-m|\le 2$. 
Indeed, by the Chapman-Kolmogorov equation, symmetry, and the fact 
$g_{nr}(j)\le 1$ for all $k,m\in \sS_{nr}$,
\[
p^{n,r}(1,k,m)\ge (nr)^{-d}\sum_{j\in \sS_{nr}}p^{n,r}(\tfrac 12, k,j+k) 
p^{n,r}(\tfrac 12, m, j+k)g_{nr}(j).
\]
Thus, by Jensen's inequality, (\ref{1.25}) yields
\[
r^dp^n(r^2,rk,rl)=p^{nr}(1,k,l)\ge \eps\qquad D\ge 1, |k-l|\le 2.
\] 
Taking $t=r^2$, this gives (\ref{1.16}). 

So we will prove (\ref{1.25}). Let $k,m\in \sS_n$ satisfy $|k-m|\le 2$ 
and set $u_t(l)=p^{n,r}(t,k,l+m)$. Define 
\[
G(t)=(nr)^{-d}\sum_{l\in \sS_{nr}}\log(u_t(l))g_{nr}(l).
\]
By Jensen's inequality, we see that $G(t)\le 0$. Further,
\[
G'(t)=(nr)^{-d}\sum_{l\in \sS_{nr}}\frac{\partial u}{\partial t}
(l)\frac{g_{nr}(l)}{u_t(l)} =-\ce^{(n),r}(u_t,\frac{g_{nr}}{u_t}).
\]
Next, note that the following elementary inequality holds (see 
page 29 of \cite{BBCK}).  
\[
\Big(\frac db -\frac ca\Big)(b-a)\le -(c\wedge d) \Big(\log \frac b{d^{1/2}}-\log \frac a{c^{1/2}}\Big)^2
+(d^{1/2}-c^{1/2})^2,\qquad a,b,c,d>0.
\] 
Applying this with $a=u_t(l),\,b=u_t(l+m),\, c=g_{nr}(l),\, d=g_{nr}(l+m)$, we have  
\begin{align*}
&G'(t)\\&=-(nr)^{2-d}
\sum_{l\in \sS_{nr}}\sum_{m\in \sS_{nr}}
\Big(\frac{g_{nr}(l+m)}{u_t(l+m)}-\frac{g_{nr}(l)}{u_t(l)}\Big)
\Big(u_t(l+m)-u_t(l)\Big)C^n(rl,r(l+m)) \\
&\ge  (nr)^{2-d}\sum_{l\in \sS_{nr}}
\sum_{m\in \sS_{nr}}(g_{nr}(l+m)\wedge g_{nr}(l))
\Big(\log \frac{u_t(l+m)}{g_{nr}(l+m)^{1/2}}-\log \frac{u_t(l)}{g_{nr}(l)^{1/2}}\Big)^2
C^n(rl,r(l+m)) \\
&\qquad -(nr)^{2-d}\sum_{l\in \sS_{nr}}
\sum_{m\in \sS_{nr}}(g_{nr}(l+m)^{1/2}-g_{nr}(l)^{1/2})^2
C^n(rl,r(l+m)) \\ 
&\ge  c(nr)^{2-d}\sum_{l\in \sS_{nr}}\sum_{j=1}^d g_{nr}(l)
\Big(\log u_t\Big(l+\frac{e^j}{nr}\Big)-\log u_t(l)+\frac 12 \Big(\Big|l_j+\frac 1{nr}\Big|-|l_j|\Big)\Big)^2 \\
&\qquad -(nr)^{2-d}\sum_{l\in \sS_{nr}}
\sum_{m\in \sS_{nr}}(g_{nr}(l+m)^{1/2}-g_{nr}(l)^{1/2})^2C^n(rl,r(l+m))=: I-II,
\end{align*}
where the last inequality is due to (A2) and the definition of $g_{nr}$
(here $e^j$ is the element of $\Z^d$ whose $k$-th 
component is $1$ if $k=j$ and $0$ otherwise). Note that 
\[
(g_{nr}(l+m)^{1/2}-g_{nr}(l)^{1/2})^2
\le c_1(|m|^2\wedge 1) (g_{nr}(l+m)+ g_{nr}(l)). 
\] 
Thus 
\begin{align*}
II \le & c_2(nr)^{2-d} \sum_{l\in \sS_{nr}}
\sum_{m\in \sS_{nr}}(g_{nr}(l+m)+ g_{nr}(l))(|m|^2\wedge 1) 
C^n(rl,r(l+m)) \\
= & 2c_2(nr)^{2-d} \sum_{l\in \sS_{nr}}\sum_{m\in \sS_{nr}} 
g_{nr}(l)(|m|^2\wedge 1)C^n(rl,r(l+m)) \\
\le & c_3\Big(\sup_{l\in \sS_{nr}}n^2\sum_{m\in \sS_{nr}}
(r^2|m|^2\wedge r^2)C^n(rl,r(l+m))\Big)\cdot (nr)^{-d} 
\sum_{l\in \sS_{nr}}g_{nr}(l) \\ 
\le & c_4\Big(\sup_{l'\in \sS_n}n^2\sum_{m'\in \sS_n}(|m'|^2\wedge 1) 
C^n(l',l'+m')\Big) \le c_5,
\end{align*}
where we used $r\le 1$ in the third inequality and (A3) in the last 
inequality. Further, 
since $(a+b)^2\geq \frac12 a^2-b^2$,
\begin{align*}
I \ge & c(nr)^{2-d}\sum_{l\in \sS_{nr}}\sum_{j=1}^d g_{nr}(l)
\Big\{\frac 12 \Big(\log u_t\Big(l+\frac{e^j}{nr}\Big)-\log u_t(l)\Big)^2
-\Big(\frac 12 \Big(\Big|l_j+\frac 1{nr}\Big|-|l_j|\Big)\Big)^2 
\Big\}\\
\ge &\frac c2 (nr)^{2-d}\sum_{l\in \sS_{nr}}\sum_{j=1}^d g_{nr}(l)
\Big(\log u_t\Big(l+\frac{e^j}{nr}\Big)-\log u_t(l)\Big)^2
-\frac{cd}4 (nr)^{-d}\sum_{l\in \sS_{nr}}g_{nr}(l)\\
\ge &\frac c2 (nr)^{2-d}\sum_{l\in \sS_{nr}}\sum_{j=1}^d g_{nr}(l)
\Big(\log u_t\Big(l+\frac{e^j}{nr}\Big)-\log u_t(l)\Big)^2-c'
\end{align*}
Combining these, we have
\begin{eqnarray*}
G'(t)&\ge &c_6(nr)^{2-d}\sum_{l\in \sS_{nr}}\sum_{j=1}^d
\Big(\log u_t\Big(l+\frac{e^j}{nr}\Big)-\log u_t(l)\Big)^2g_{nr}(l)-c_5 \\
&\ge & c_7(nr)^{-d}\sum_{l\in \sS_{nr}}(\log u_t(l)-G(t))^2g_{nr}(l)-c_5,
\end{eqnarray*}
where we used Lemma \ref{1.19} in the last inequality. 
Given these estimates, the rest of the proof is very similar to that
of \cite[Lemma 4.4]{BaKu08}. 
\qed

\begin{rem}
{\rm 
There is an error in the proof of \cite[Lemma 4.4]{BaKu08}. 
The estimate $|g_D(l+e)-g_D(l)|\le c_1D^{-1}|e| (g_D(l+e)\wedge g_D(l))$ 
in page $2051$, line 23, is not true when $D\ll |e|$. However, 
one can easily fix the proof by arguing as in the proof here.
}
\end{rem}

The next lemma can be proved exactly in the same way as \cite[Lemma 4.5]{BaKu08}. 

\begin{lem}\label{offdiag}
Given $\delta>0$ there exists $\kappa$ such that for each $n\in \N$, 
if $x,y\in \sS_n$ and $C\subset \sS_n$ with $\mbox{\rm dist}\,(x,C)$ and 
$\mbox{\rm dist}\,(y,C)$ both larger than $\kappa t^{1/2}$ where 
$t\in (n^{-1},1]$, then 
\[ 
\P^x(Y^{(n)}_t=y, T^n_C\leq t)\leq \delta t^{-d/2}n^{-d}.
\]
\end{lem}

\noindent 
{\sc Proof of Proposition \ref{szhol}.} 
We have from Lemma \ref{1.20} that there exists $\varepsilon$ such that
\[ \P^x(Y^{(n)}_t=y)=p^n(t,x,y)\mu_y^n\geq \varepsilon t^{-d/2}n^{-d} \]
if $|x-y|\leq 2t^{1/2}$. If we take $\delta=\varepsilon/2$ in Lemma 
\ref{offdiag}, then provided $r>(\kappa +1) t^{1/2}$, we have
\[ 
\P^x(Y^{(n)}_t=y, \tau^n_{B_n(x_0,r)}\leq t)\leq \frac{\eps}{2} 
t^{-d/2}n^{-d}. 
\] 
Subtracting,
\[ 
\P^x(Y^{(n)}_t=y, \tau^n_{B_n(x_0,r)}>t)\geq \frac{\eps}{2} t^{-d/2}n^{-d} 
\]
if $|x-y|\leq t^{1/2}$, which is equivalent to what we want.
\qed

For $(t, x)\in [0, 1]\times \sS_n$
and $r>0$ let $Q^n(t, x, r):=[t, \, t+\gamma r^2] \times B_n(x, r)$,
where $\gamma:=\gamma(1/2, 1/2)<1$.  Here $\gamma(1/2, 1/2)$ is the 
constant in (\ref{eq:tighty}) corresponding to $A=B=1/2$. 

Given the above estimates, we can prove the uniform H\"older continuity of 
the heat kernel $p^n (t, x, y)$ similarly to \cite[Theorem 4.9]{BaKu08}. 

\begin{thm}\label{4.15}
There are constants $c>0$ and $\beta>0$ (independent of $R,n$) such that for
every $0<R\le 1 $, every $n\ge 1$, and every bounded parabolic function
$q$ in $Q^n(0, x_0, 4R)$, 
\begin{equation}\label{eqn:holder1}
|q(s, x) -q(t, y)| \leq c \,  \| q \|_{\infty, R} \, R^{-\beta} \, 
\left( |t-s|^{1/2} + |x-y| \right)^\beta
\end{equation}
holds for $(s, x), \, (t, y)\in Q^n(0, x_0, R)$, where $\| q \|_{\infty, R}
:=\sup_{(t,y)\in [0, \, \gamma (4R)^2 ] \times \sS_n } |q(t,y)|$. In 
particular, for the transition density function $p^n (t, x, y)$ of $Y^{(n)}$, 
\begin{equation}\label{eqn:holder2}
 |p^n (s, x_1, y_1) -p^n (t, x_2, y_2)| \leq c \, t_0^{-(d+\beta)/2}
\left( |t-s|^{1/2} + |x_1-x_2|+|y_1-y_2| \right)^\beta,
\end{equation}
for any $n^{-1}<t_0<1$, $t, \, s \in [t_0, \, 1]$ and $(x_i, y_i)\in
\sS_n\times \sS_n$ with $i=1, 2$. 
\end{thm}

\proof Given the above estimates, we can prove 
the analogues of 
Corollary 4.6 and Lemma 4.7 in \cite{BaKu08} exactly in the same way
as is done there. 
Thus the proof of Theorem \ref{4.15} is almost the same as that of 
\cite[Theorem 4.9]{BaKu08} except for the following small change.

The following computation is needed to obtain 
the first inequality of (4.13) in \cite{BaKu08}:  
\[
\sup_{z\in B_n(x,r)}n^2\sum_{y\in\sS_n\setminus \overline{B_n(x, s)}}
C^n(z,y)\le \Big(\frac s2\Big)^{-2} \sup_{z\in B_n(x,r)}\sum_{y\in\sS_n}
\Big(|z-y|^2\wedge 1\Big)C^n(z,y)n^2\le \frac{c_2}{s^2}
\]
where (A3) is used in the last inequality (note that $2r\le s\le 1$).  \qed

\section{Weak convergence of the process}\label{S:wc}
Recall that $Y^{(n)}_t$ are the continuous time Markov chains on $\sS_n$ 
corresponding to $(\ce^n,\cf^n)$ in (\ref{DFe_n}) and (\ref{DFf_n}).  
Since the state space of $Y^{(n)}$ is $\sS_n$ while the limit process
will have $\R^d$ as its state space, we need to exercise some care with  
the domains of the functions we deal with. First, if $g$ is defined on 
$\R^d$, we define $R_n(g)$ to be the restriction of $g$ to $\sS_n$:
\[ 
R_n(g)(x)= g(x), \qquad x\in \sS_n. 
\]
If $g$ is defined on $\sS_n$, we define $E_ng$ to be the  extension of $g$ to 
$\R^d$  defined by
$$
E_ng(x)=g([x]_n), 
$$ 
where $[x]_n=([nx_1]/n,[nx_2]/n,\ldots,[nx_d]/n)$ for $x=(x_1,x_2,\ldots,x_d) 
\in \rd$. 

In order to consider the convergence of the processes and to identify 
the limit process, we need to show the  convergence of the semigroups of the 
Dirichlet forms $(\form^n,\dom^n)$ in an appropriate sense.  
To this end, we now prepare some notation to specify a 
condition under which the convergence holds. 
For $n\in {\mathbb N}$, set 
$$
|x-y|_n:=n|x_1-y_1|+n|x_2-y_2|+\cdots+n|x_d-y_d| \ (\in \N), \quad 
{\rm for} \ x,y \in \Sn.
$$
Note that $1\le |x-y|_n \le d n|x-y|$ holds for any $x,y\in \Sn$ with 
$x\not=y$, where $|x-y|$ is the Euclidean distance between $x$ and $y$.

{\it A shortest path $\sigma$ from $x$ to $y$} is a sequence of points $p_i \in \Sn$ 
for $i=0, 1,2,\ldots, k=|x-y|_n$, which we denote by $\sigma=\sigma(p_0, \ldots, p_k)$, 
so that $p_0=x, p_k=y$ and for any $\ell=0, 1,\ldots, k-1$, there exists 
$j\in \{1,2,\ldots, 2d\}$ such that 
$$
p_{\ell}=p_{\ell+1} +\frac 1n \mbox{\boldmath  $\alpha$}_j, 
$$
where $\mbox{\boldmath  $\alpha$}_i=\mbox{\boldmath  $e$}_i$ if $i=1,2,\ldots, 
d$ and $\mbox{\boldmath  $\alpha$}_i=-\mbox{\boldmath  $e$}_{i-d}$ if 
$i=d+1,\ldots, 2d$. Let ${\cal P}(x,y)$ be the set of all shortest paths $\sigma$ 
from $x$ to $y$. The number of all such shortest paths $\sigma$ is 
$$
\Pi(x,y):=\frac{(|x-y|_n) !}{ \bigl(n|x_1-y_1|\bigr) ! \bigl(n|x_2-y_2|\bigr) ! 
\cdots \bigl(n|x_d-y_d|\bigr) !}. 
$$
For $\sigma \in {\cal P}(x,y)$, define a function $D_{\sigma}$ 
defined on $\Sn\times\Sn$ as follows:
$$
D_{\sigma}(w,z) :=
\begin{cases}
1, & \text{ if there exists  $\ell$  such that  
$w=p_{\ell}$  and $ z=p_{\ell+1}$}, \\ 
0, &\text{ otherwise}. 
\end{cases} 
$$
For any function $u$ defined on $\Sn$ and for any $x,y\in \Sn$, we easily 
see that 
$$
u(x)-u(y)= \dis{
\frac1{\Pi(x,y)} \sum_{\sigma \in {\cal P}(x,y)} \sum_{z,w\in \Sn} 
D_{\sigma}(w,z) \bigl(u(w)-u(z)\bigr). } 
$$
Now let 
$$
P^{x,y}(w,z)=\frac 1{\Pi(x,y)} \sum_{\sigma \in {\cal P}(x,y)} 
D_{\sigma}(w,z).
$$
For $h\in {\mathbb R}$, $x\in \rd$  and $i=1,2,\ldots,d$, let 
$$
\nabla_h^i u(x) =\frac{u(x+h\mbox{\boldmath $e$}_i)-u(x)}{h}. 
$$ 
We then have the following.
\begin{lem}\label{expr_pxy}
$$ \dis{
u(x)-u(y)= \frac 1n \sum_{i=1}^d \sum_{z\in \Sn}
\Bigl(P^{x,y}(z+\mbox{\boldmath $e$}_i/n, z)-
P^{x,y}(z, z+\mbox{\boldmath $e$}_i/n)\Bigr)\nabla_{1/n}^i u(z)}.
$$
\end{lem}
\proof We have 
\begin{align*}
\dis{\sum_{w\in \Sn} D_{\sigma}(w,z)}& \bigl(u(w)-u(z) \bigr)\\
& \dis{ = \sum_{i=1}^{2d} D_{\sigma}(z+
\mbox{\boldmath $\alpha$}_i/n, z) \bigl(u(z+ 
\mbox{\boldmath $\alpha$}_i/n)-u(z) \bigr)  } \\
& \dis{ = \sum_{i=1}^{d} \Bigl\{ D_{\sigma}(z+
\mbox{\boldmath $e$}_i/n, z) \bigl(u(z+
\mbox{\boldmath $\alpha$}_i/n) - u(z) \bigr) } \\
&  \quad \quad \quad  \dis{ 
+D_{\sigma}(z -
\mbox{\boldmath $e$}_i/n, z)  \bigl(u(z-
\mbox{\boldmath $e$}_i/n) -u(z)\bigr) \Bigr\} } \\
&= \dis{ \frac 1n
\sum_{i=1}^{d} \Bigl\{ D_{\sigma}(z+
\mbox{\boldmath $e$}_i/n, z)  \nabla^i_{1/n} u(z) 
 -D_{\sigma}(z -
\mbox{\boldmath $e$}_i/n, z)  \nabla_{-1/n}^i u(z)  \Bigr\}. } \\
\end{align*}
So
\begin{align*} u(x)&-u(y)\\ 
&= \dis{\sum_{z\in \Sn} \frac 1{\Pi(x,y)} \sum_{\sigma \in {\cal P}(x,y)} 
\sum_{w\in \Sn} D_{\sigma}(w,z) \bigl(u(w)-u(z)\bigr) } \\
&= \dis{ \frac 1n \sum_{i=1}^d \sum_{z\in \Sn}
\frac 1{\Pi(x,y)} \sum_{\sigma \in {\cal P}(x,y)} 
\Bigl( D_{\sigma}(z+
\mbox{\boldmath $e$}_i/n, z)  \nabla^i_{1/n} u(z)  
-D_{\sigma}(z - \mbox{\boldmath $e$}_i/n, z)  
\nabla_{-1/n}^i u(z)  \Bigr) } \\
&= \dis{ \frac 1n \sum_{i=1}^d \sum_{z\in \Sn}
\Bigl(P^{x,y}(z+\mbox{\boldmath $e$}_i/n, z) \nabla_{1/n}^i u(z) 
-P^{x,y}(z-\mbox{\boldmath $e$}_i/n, z) \nabla_{-1/n}^i u(z) \Bigr). } \\
\end{align*}
Moreover, for each $i=1,2,\ldots,d$, and $x,y\in \Sn$, 
$$
\begin{array}{lcl}
\dis{ \sum_{z\in \Sn}
P^{x,y}(z-\mbox{\boldmath $e$}_i/n, z) \nabla_{-1/n}^i u(z)} 
&=& \dis{ \sum_{z\in \Sn}
P^{x,y}(z, z+\mbox{\boldmath $e$}_i/n) \nabla_{-1/n}^i u(z+
\mbox{\boldmath $e$}_i/n)} \\
& &\vspace{-6pt} \\
&=& -n \dis{ \sum_{z\in \Sn} P^{x,y}(z, z+\mbox{\boldmath $e$}_i/n) 
\bigl(u(z)-u(z+\mbox{\boldmath $e$}_i/n)\bigr)} \\
& &\vspace{-6pt} \\
&=& \dis{ \sum_{z\in \Sn} P^{x,y}(z, z+\mbox{\boldmath $e$}_i/n) 
\nabla_{1/n}^i u(z). } \\
\end{array}
$$
We thus obtain the desired equality. \qed  
\begin{rem}
{\rm 
Here $P^{x,y}(\cdot,\cdot)$ is defined by averaging over the 
set of  all shortest paths between $x$ and $y$. However, we could take an
average over other collections of paths. 
Let $x=(x_1,\cdots, x_d), y=(y_1,\cdots, y_d)$. 
Other possible collections of paths  are the following:\\
(i) Let $H(x,y)$ be the $d$-dimensional cube whose vertices consist of  
$\{(z_1,\cdots,z_d): z_i $ is either $x_i$ or $y_i$ for $i=1,\cdots, d\}$. 
Let ${\cal P}(x,y)$ be the set of shortest paths between $x$ and $y$ 
that consist of a union of the edges of $H(x,y)$, and 
take the average over  ${\cal P}(x,y)$. 
In this case $\Pi(x,y)$ in the definition of  $P^{x,y}(\cdot,\cdot)$ is $d\,!$.\\
(ii) Let $L_{x,y}$ be the union of the line segment from $x$ to $(y_1,x_2,\cdots, x_d)$, 
the line segment from $(y_1,x_2,\cdots, x_d)$ to $(y_1,y_2, x_3\cdots, x_d)$, $\cdots$, 
and the  line segment from $(y_1,\cdots, y_{d-1},x_d)$ to $y$. 
Set ${\cal P}(x,y)=\{L_{xy}\}$ and $\Pi(x,y)=1$. This was used in \cite{BaKu08}.
}
\end{rem}

Next, 
let us fix a decreasing sequence $\{\vareps_n\}$ such that 
$1\ge \vareps_n\searrow 0$, and define functions $C_C^n(x,y)$, 
$C_J^n(x,y)$ on $\sS_n\times\sS_n$ as follows:
$$
C_C^n(x,y)  :=
\left\{ 
\begin{array}{ll} 
C^n(x,y), &  {\rm if} \quad |x-y|\le \vareps_n, \\
0,        & {\rm otherwise}, \\
\end{array} 
\right. 
$$
and $C_J^n(x,y) := C^n(x,y)-C_C^n(x,y), \quad x,y \in \sS_n$.

\bigskip
Now define the following Dirichlet forms 
corresponding to the conductances $C_C^n(x,y)$ and $C_J^n(x,y)$, 
which we consider as the `continuous part' and the `jump part' of the Dirichlet form 
$(\form^n,\dom^n)$; for $f\in L^2(\sS_n,\mu^n)$, 
$$
\left\{
\begin{array}{lcl}
\form_C^n(f,g)&:=& \dis{\frac{n^{2-d}}2\sum_{x,y \in \sS_n} 
\bigl(f(x)-f(y)\bigr) \bigl(g(x)-g(y)\bigr) C_C^n(x,y),} \\ 
& & \vspace{-6pt} \\
\form_J^n(f,g)&:=& \dis{\frac{n^{2-d}}2\sum_{x,y \in \sS_n} 
\bigl(f(x)-f(y)\bigr) \bigl(g(x)-g(y)\bigr) C_J^n(x,y).} \\ 
\end{array}
\right. 
$$
Then clearly $\form^n(f,g)=\form_C^n(f,g)+\form_J^n(f,g)$. 

Using Lemma \ref{expr_pxy}, we can write down $\form_C^n(u,v)$ as follows:
\begin{equation}
\begin{array}{lcl} \form_C^n(u,v) 
&=& \dis{ \frac{n^{2-d}}2 \sum_{x,y\in \Sn} 
\bigl(u(x)-u(y)\bigr)\bigl(v(x)-v(y)\bigr) C^n_C(x,y) } \\

& & \vspace{-6pt} \\
&=& \dis{ \frac 1{2n^{d}} \sum_{x,y\in \Sn} \sum_{i,j=1}^d \sum_{z,w\in \Sn} 
 \Bigl(P^{x,y}(z+\mbox{\boldmath $e$}_i/n, z)-
P^{x,y}(z, z+\mbox{\boldmath $e$}_i/n)\Bigr) } \\

& & \vspace{-6pt} \\
& &  \quad \times \dis{ \Bigl(P^{x,y}(w+\mbox{\boldmath $e$}_j/n, w)-
P^{x,y}(w, w+\mbox{\boldmath $e$}_j/n)\Bigr)\nabla_{1/n}^i u(z) 
\nabla_{1/n}^j v(w) C^n_C(x,y).  }
\end{array}
\label{eta-form}
\end{equation}
For $i,j=1,2,\ldots,d$ and $w,z \in \Sn$, set
$$
\begin{array}{lcl} {G}^n_{ij}(w,z)
&:=&  \dis{ \sum_{x,y\in \Sn} \Bigl(P^{x,y}(z+\mbox{\boldmath $e$}_i/n, z)-
P^{x,y}(z, z+\mbox{\boldmath $e$}_i/n)\Bigr) } \\

& & \vspace{-6pt} \\
& & \quad \quad \times \dis{ 
\Bigl(P^{x,y}(w+\mbox{\boldmath $e$}_j/n, w)-
P^{x,y}(w, w+\mbox{\boldmath $e$}_j/n)\Bigr) C^n_C(x,y); } \\
\end{array}
$$
then we see that  
\begin{equation}
\form_C^n(u,v)=\frac 1{2n^{d}} \sum_{i,j=1}^d \sum_{w,z \in \Sn} \nabla^i_{1/n}u(z)
\nabla_{1/n}^j v(w) {G}^n_{ij}(w,z).
\label{(6)}
\end{equation}
Let 
\begin{equation}
F_{ij}^n(z)=\sum_{w\in \Sn} {G}_{ij}^n(w,z), 
\quad z\in \Sn, \ i,j=1,2,\ldots,d. 
\label{(7)}
\end{equation}
Note that if (A4) below holds, then by 
the fact that $C_C^n(x,y)=0$ for $|x-y| 
> \vareps_n$, we have $F_{ij}^n \in L^1(\Sn,\mu^n)$.

\bigskip

From now on, we extend the conductances $C^n(x,y)$ to $\rd\times\rd$ as 
follows:
$$
C^n(x,y)=C^n([x]_n,[y]_n) \quad {\rm for} \ x,y \in \rd.
$$
We extend $C^n_C( \cdot,\cdot), C^n_J( \cdot,\cdot)$ to $\rd\times \rd$ and 
extend $F_{ij}^n(\cdot)$ to $\rd$ similarly. 

\bigskip
We now give an assumption needed to obtain weak convergence of the processes. 
\begin{itemize}
\item[(A4)] {\it There exist a decreasing sequence $\{\vareps_n\}$ satisfying 
$1/n \le \vareps_n \le 1$ and $\vareps_n \searrow 0$, symmetric matrix-valued 
functions $a(x)=(a_{ij}(x))$ on $\rd$, and symmetric functions $j(x,y)$ on 
$\rd\times\rd\backslash D$ so that for any $i,j=1,2,\ldots,d$, the functions 
$F_{ij}^n(x)$ are uniformly bounded and converge to $a_{ij}(x)$ 
locally in $L^1(\rd)$, and 
$$
\lambda^{-1}|\xi|^2 \le \sum_{i,j=1}^d \xi_i\xi_j a_{ij}(x) \le 
\lambda |\xi|^2, \quad x,\xi \in \rd, 
$$
for some $\lambda>0$. Further, for each $N>1$, 
the measures 
\begin{equation} \label{WCA3}
n^{d+2} C^n(x,y){\bf 1}_{[N^{-1},N]}(|x-y|)dx\,dy 
\longrightarrow j(x,y){\bf 1}_{[N^{-1},N]}(|x-y|)dx\,dy
\end{equation}
{ weakly} as $n\to\infty$.
}
\end{itemize}

\begin{rem}\label{RR5.3}{\rm 
Here \eqref{WCA3} refers to the weak convergence of the measures on the left
to the measures on the right. Saying that the $F^n_{ij}$ are uniformly bounded
and converge locally in $L^1$ means that $\sup_{i,j,n}||F^n_{ij}||_\infty<\infty$
and for every compact set $B$, 
$$\int_B |F^n_{ij}(x)-a_{ij}(x)|\, dx\to 0.$$
Since the $F^n_{ij}$ are uniformly bounded, the convergence locally in $L^1$
is equivalent to the convergence in measure on each compact set. In particular,
a subsequence will converge almost 
everywhere.  
}
\end{rem}

From (A3) and (A4), we have 
$$
\sup_x \int_{y\not=x} \hspace{-5pt} \bigl(1\wedge|x-y|^2\bigr)
j(x,y) dy \le \int_{y\not=x} \hspace{-5pt} \bigl(1\wedge|x-y|^2\bigr)
\varphi(|x-y|)dy =\int_{h\not=0} \hspace{-5pt} 
\bigl(1\wedge|h|^2\bigr) \varphi(|h|)dh <\infty.
$$
Since $a$ is uniformly elliptic, if we define 
$$
\begin{array}{lcl}
\form(f,g) &:=& \form_C(f,g) + \form_J(f,g) \\

& & \vspace{-6pt} \\
 &:=& \dis{\frac 12 \int_{\rd} \nabla f(x) \cdot a(x)\nabla g(x) dx
+ \frac 12 \iint_{x\not=y} \bigl(f(x)-f(y)\bigr)\bigl(g(x)-g(y)\bigr) j(x,y)dx\,dy}, \\
\end{array}
$$
then $(\form,C_c^1(\rd))$ is a closable Markovian form on $L^2(\rd, dx)$. 
Denote the closure by $(\form,\dom)$.
\begin{lem}\label{thm:domlelem} 
Let $W^{1,2}(\rd):=\{f\in L^2(\rd, dx):  \nabla f \in L^2(\rd, dx)\}$. Then,
\begin{equation}\label{w12maxd}
\{f\in L^2(\rd, dx): \form(f,f)<\infty\}=W^{1,2}(\rd) =\sF.
\end{equation}
Further, if 
$(\form,\dom')$ is a regular Dirichlet form on $L^2(\rd, dx)$, 
then $\dom'=W^{1,2}(\rd)$. 
\end{lem}
\proof  Let $f\in L^2$ be such that 
$\form (f,f)<\infty$. Then, $\form_C(f,f)<\infty$ and  
$\form_C(f,f)$ is comparable to 
$\|\nabla f\|_2^2$,
so $f\in W^{1,2}(\rd)$. On the other hand, 
suppose $f\in W^{1,2}(\rd)$. Then, $\sE_J(f,f)\le \sE_\varphi (f,f)$, where $\varphi$ is given in (A3) and $\sE_\varphi$ is the Dirichlet form for the symmetric L\'evy process with L\'evy measure $\varphi (|h|)dh$.  
By  
the L\'evy-Khintchine formula (see {\it e.g.} (1.4.21) in \cite{FOT}), 
the characteristic function $\psi$ of the process is given by 
$$
\psi(u)=\int_{\rd} \Bigl(1-\cos \bigl(u\cdot h\bigr)\Bigr) \varphi(|h|) dh, 
\quad u\in \rd.
$$
According to (A3), 
we have,  
\begin{align*}
\psi(u) &=\int [1-\cos(u\cdot h)]\varphi(|h|)\, dh\\
 & \leq c_1 \int [|u|^2|h|^2\land 1] \varphi(|h|)\, dh\\
 &\leq c_2(|u|^2 + 1). \\
 \end{align*}
Using Plancherel's theorem, for $f\in C^2_c(\R^d)$,

\begin{align*}
\sE_\varphi(f,f) & = \frac 12 \iint_{y\not=x} \bigl(f(x+h)-f(x)\bigr)^2 
\varphi(|h|)dhdx \\
 & =\int \bigl| \widehat{f}(u)\bigr|^2 \psi(u)du \\
&\leq c_2\int (1+|u|^2) |\widehat f(u)|^2\, 
du=c_3(\|f\|_2^2 +\|\nabla f\|_2^2).
\end{align*}
Here $\hat f$ is the  Fourier transform of $f$. 
A limit argument shows that
\begin{equation}\label{plause2}
\sE_\varphi(f,f) \le c_4(\|f\|^2_2+\|\nabla f\|_2^2) 
\end{equation}
for $f\in W^{1,2}(\R^d)$. 
Since $\form_C(f,f)$ is comparable to 
$\|\nabla f\|_2^2$,
adding shows that $\form(f,f)<\infty$,  
and the first equality in (\ref{w12maxd}) is proved. 
Now suppose $(\form,\dom')$ is a regular Dirichlet form on $L^2(\rd, dx)$;  
then since $W^{1,2}(\R^d)$ is the maximal domain (due to the first equality 
in (\ref{w12maxd})), we have $\sF'\subset W^{1,2}(\R^d)$. 
From the above results, we know that the $(\sE(\cdot,\cdot)+
\|\cdot\|_2^2)$-norm is comparable to the $W^{1,2}$-norm on $W^{1,2}(\R^d)$. 
Using this, 
we see that $(\|\nabla \cdot\|_2^2,\dom')$ is a regular Dirichlet form. This 
implies $\sF'= W^{1,2}(\R^d)$ (so $W^{1,2}(\rd) =\sF$ as well) and  
the proof is complete.   
\qed 

\medskip

Under the above set-up we have the following, which is the  main 
theorem of this paper. 
\begin{thm}\label{clt-main}
Suppose (A1)-(A4) hold. Then for each $x$ and each $t_0$ the 
$\P^{[x]_n}$-laws of $\{Y^{(n)}_t; 0\leq t\leq t_0\}$ converge weakly 
with respect to the topology of the space $D([0,t_0], \R^d)$. 
If $Z_t$ is the canonical process on $D([0, t_0), \R^d)$ and $\P^x$ 
is the weak limit of the $\P^{[x]_n}$-laws of $Y^{(n)}$, then the process 
$\{Z_t, \P^x\}$ is the symmetric Markov process corresponding to the Dirichlet 
form $\ce$ with domain $W^{1,2}(\R^d)$. 
\end{thm}

\section{Proof of Theorem \ref{clt-main}}\label{S:proof}

In this section, we will prove Theorem \ref{clt-main}. We first extend 
$\sE^n$ and define a quadratic form on $L^2(\rd, dx)$.
Define 
$$
{\cal H}_n:=\Bigl\{E_n u:    \,u \mbox{ is a function on }\sS_n\Bigr\}
\cap L^2(\rd, dx).
$$
For $f=E_n u\in {\cal H}_n$, define  
$$
\tilde{\form}^n(f,f) = \dis{
\frac{n^{2+d}}2 \iint_{x\not=y} \bigl(f(x)-f(y)\bigr)^2  C^n(x,y) dx\, dy. } 
$$
Then we see 
\begin{equation}
\begin{array}{lcl}
\tilde{\form}^n(f,f)&=& 
\dis{ \frac{n^{2+d}}2 \sum_{w_1,w_2\in \Sn} \bigl(u(w_1)-
u(w_2)\bigr)^2 C^n(w_1,w_2) \bigl(n^{-d}\bigr)^2 } \\

& & \vspace{-6pt} \\
&=& \dis{ \frac{n^{2-d}}2 \sum_{w_1,w_2\in \Sn} \bigl(u(w_1)-
u(w_2)\bigr)^2 C^n(w_1,w_2) =\form^n(u,u).}  \\
\end{array}
\label{eqn:form-equiv}
\end{equation}

\bigskip

\noindent {\sc Proof of Theorem \ref{clt-main}}. Let $U_n^\lam$ be the 
$\lam$-resolvent for $Y^{(n)}$; this means that 
\[ 
U^\lam_n h(x)=\E^x\int_0^\infty e^{-\lam t} h(Y^{(n)}_t)\, dt 
\]
for $x\in \sS_n$ and $h: \sS_n\to\R$. The first step is to show 
that any subsequence $\{n_j\}$ has a further 
subsequence $\{n_{j_k}\}$  such that $U^\lam_{n_{j_k}}(R_{n_{j_k}}f)$ 
converges uniformly on compacts whenever $f\in C_c(\mathbb R^d)$, 
that is, $f$ is continuous with compact support. Given Proposition 
\ref{thm:4.2} and Theorem \ref{4.15}, the proof of this is very similar to 
that of \cite[Proposition 6.2]{BaKu08}, and we refer the reader to that paper. 

Now suppose we have a subsequence $\{n'\}$ such that the 
$U^\lam_{n'} (R_{n'}f)$ are equicontinuous and converge uniformly on 
compacts whenever $f\in C_c(\mathbb R^d)$. Fix such an $f$ and let $H$ be 
the limit of $U^\lam_{n'}(R_{n'}f)$. Let $g\in C^2_c(\R^d)$ and
write
$\angel{f,g}:=\int_{\mathbb R^d}
f(x)g(x)dx$.

In the following, we drop the primes for legibility. 
Set $u_n=U^\lam_n (R_nf)$ for $\lam>0$. 
We will  
prove that 
\begin{equation}\label{WCElqpw}
H\in W^{1,2}(\rd)\qquad \mbox{and}\qquad  
\sE^n(u_n,g)\to \sE(H,g)
\end{equation}
along some subsequence.
Once we have \eqref{WCElqpw}, then
\begin{align*}
\sE(H,g)&=\lim \sE^n(u_n,g)=\lim
(\angel{f,g}_n-\lam \angel{u_n,g}_n)\\
&=\angel{f,g}-\lam \angel{H,g},
\end{align*}
the limit being taken along the subsequence and
where $\angel{h_1,h_2}_n=n^{-d}\sum_{x\in \sS_n}h_1(x)h_2(x)$ 
for $h_1,h_2: \sS_n\to \mathbb R$. 
By \eqref{WCElqpw}, $H\in W^{1,2}(\R^d)$, and the equality
\begin{equation}\label{WCERB1}
\sE(H,g)=\angel{f,g}-\lam\angel{H,g}
\end{equation}
holds for all $g\in C^2_c(\R^d)$. By
Lemma \ref{thm:domlelem}, $C^2_c(\R^d)$ is dense in $W^{1,2}(\R^d)$ with
respect to the norm $(\sE(\cdot, \cdot)+\|\cdot\|_2^2)^{1/2}$, and
so \eqref{WCERB1} holds for all $g\in W^{1,2}(\R^d)$. 
Since $W^{1,2}(\R^d)$ is the maximal domain due to 
(\ref{w12maxd}), 
this implies that $H$ is the $\lambda$-resolvent of $f$ for the process 
corresponding to $(\sE, W^{1,2}(\R^d))$, that is, $H=U^\lam f$.
We can then conclude that the full sequence $U^\lam_n (R_{n}f)$ 
(without the primes) converges 
to $U^\lam f$ whenever $f\in C_c(\mathbb R^d)$. The assertions about the 
convergence of $\P^{[x]_n}$ then follow as in \cite[Proposition 6.2]{BaKu08}.
The rest of the proof will be devoted to proving (\ref{WCElqpw}). 

\medskip
\noindent {\it The jump part.}   

\medskip

This part of the proof is similar to that 
of \cite[Theorem 4.1]{BKK}. 
We know
\begin{equation}\label{WCE3}
\sE^n(u_n, u_n)=\angel{R_nf,u_n}_n-\lam \norm{u_n}_{2,n}^2.
\end{equation} 
Since $||\lambda u_n||_{2,n}^2
=||\lambda U^{\lambda}_nR_n f||_{2,n} \le ||R_nf||_{2,n}\le 
\sup_n ||R_nf||_{2,n}$ (note that $\sup_n ||R_nf||_{2,n}<\infty$ 
because $\lim_{n\to\infty}||R_nf||_{2,n}=\|f\|_2$ for $f\in C_c(\rd)$), 
the right hand side of  
(\ref{WCE3}) is bounded by 
$$
\bigl|\angel{R_nf,u_n}_n \bigr| + 
\lambda || u_n||_{2,n}^2  
\le \dis{ \frac{1}{\lambda} ||R_nf||_{2,n} 
||\lambda u_n||_{2,n} + 
\frac{1}{\lambda}||\lambda u_n||_{2,n}^2  } 
\le \dis{\frac{2}{\lambda}\sup_n ||R_nf||_{2,n}^2.} 
$$ 
This tells us that $\{\form^n(u_n,u_n)\}_n$ is uniformly bounded. 

Since the $u_n$ are equicontinuous and converge uniformly to $H$ on 
$\ol{B(0, N)}$ for $N>0$, using (\ref{WCA3}), we have 
\begin{align*}
\int\int_{N^{-1}\le |y-x|\le N} &(H(y)-H(x))^2 j(x,y)\, dy\, dx\\
&\leq \limsup_{n\to \infty} n^{2-d}\sum_{{x,y\in \sS_n}
\atop{N^{-1}<|y-x|\le N}}(u_n(y)-u_n(x))^2 C^n(x,y)\\
&\leq \limsup_n \sE^n(u_n, u_n)\leq c<\infty.
\end{align*}
Letting $N\to \infty$, we have
\begin{equation}\label{WCE2}
\sE_J(H,H)<\infty.
\end{equation}

Fix a function $g$ on $\sS_n$  
with compact support and choose $M$ large enough so that
the support of $g$ is contained in $B(0,M)$. Then
\begin{align*}
&\Big| n^{2-d}\sum_{{x,y\in \sS_n}\atop{|y-x|> N}} 
(u_n(y)-u_n(x))(g(y)-g(x))C^n(x,y) \Big|\\
&\leq \Big(n^{2-d}\sum_{x,y\in \sS_n} (u_n(y)-u_n(x))^2 C^n(x,y)\Big)^{1/2} 
 \Big(n^{2-d}\sum_{{x,y\in \sS_n}\atop{|y-x|> N}} 
(g(y)-g(x))^2 C^n(x,y)\Big)^{1/2}.
\end{align*}
The first factor is $(\sE^n(u_n, u_n))^{1/2}$, while the second factor is 
bounded by 
$$
2\norm{g}_\infty \Big(n^{2-d}\sum_{x\in B(0,M)\cap \sS_n}
\sum_{|y-x|> N} C^n(x,y)\Big)^{1/2},
$$
which, in view of \eqref{cond-bound}, will be small if $N$ is large. 
Similarly,
\begin{align*}
&\Big| n^{2-d}\sum_{{x,y\in \sS_n}\atop{|y-x|< N^{-1}}} 
(u_n(y)-u_n(x))(g(y)-g(x))C^n(x,y) \Big|\\
&\leq \Big(n^{2-d}\sum_{x,y\in \sS_n} (u_n(y)-u_n(x))^2 C^n(x,y)\Big)^{1/2} 
\cdot \Big(n^{2-d}\sum_{{x,y\in \sS_n}\atop{|y-x|< N^{-1}}} 
(g(y)-g(x))^2 C^n(x,y)\Big)^{1/2}.
\end{align*}
The first factor is as before, while the second is bounded by
$$\norm{\nabla g}_\infty 
\Big(n^{2-d}\sum_{x\in B(0,M)\cap \sS_n}
\sum_{|y-x|< N^{-1}} |y-x|^2C^n(x,y)\Big)^{1/2}.$$
In view of \eqref{cond-bound}, the second factor will be small if $N$ is large.

Using \eqref{WCE2}, we have that
$$
\Big|\int\int_{|y-x|\notin [N^{-1}, N]}
(H(y)-H(x))(g(y)-g(x)) j(x,y)\, dy\, dx\Big|
$$
will be small if $N$ is taken large enough. 

By \eqref{WCA3} and the fact that the $U^\lam_nf$ are equicontinuous and 
converge to $H$ uniformly on compacts, if we take $n$ large enough so that
$\vareps_n\le N^{-1}$, we have 
\begin{align*}
& \hspace{-20pt} n^{2-d}\sum_{{x,y\in \sS_n}\atop{N^{-1}\le |y-x|\le N}} 
 (u_n(y)-u_n(x))(g(y)-g(x))C^n_J(x,y) \\
=&n^{2-d}\sum_{{x,y\in \sS_n}\atop{N^{-1}\le |y-x|\le N}} 
(u_n(y)-u_n(x))(g(y)-g(x))C^n(x,y)\\ 
\to & \int\int_{N^{-1}\le |y-x|\le N} (H(y)-H(x))(g(y)-g(x))j(x,y)\, dy\, dx.
\end{align*}
It follows that
\begin{equation}\label{WCE4}
\sE_J^n(u_n,g)\to \sE_J(H,g),
\end{equation}
which takes care of  the jump part of (\ref{WCElqpw}). 

\medskip
\noindent{\it The continuous part.}             

\medskip
{\it Step 1.} First we show that $H\in W^{1,2}(\rd)$. 

As in the discussion of the jump part, we know $\{\form^n(u_n,u_n)\}_n$ is 
uniformly bounded.  
On the other hand, making use of the assumption (A2), we see 
$$
\tilde{\form}^n(E_nu_n,E_nu_n) = \form^n(u_n,u_n) 
\ge c \form_{NN}^n(u_n,u_n) = c \tilde{\form}_{NN}^n(E_nu_n,E_nu_n). 
$$
Therefore, for $f\in C_c^1(\rd)$, the sequence 
$\{\tilde{\form}_{NN}^n(E_nu_n,E_nu_n)\}_n$ is uniformly bounded with 
respect to $n$. Letting $Q_n(w)= \prod_{i=1}^d [w_i,w_i+1/n)$, 
we see that for any $i=1,2,\ldots,d$, 
$$
\begin{array}{lcl}
\tilde{\form}^n_{NN}(E_nu_n,E_nu_n) 
&=& \dis{\frac{n^{2+d}}{2} \iint_{x\not=y} 
\bigl(E_nu_n(x)-E_nu_n(y)\bigr)^2 C_{NN}^n(x,y)dxdy } \\
& & \vspace{-7pt} \\
&=& \dis{ \frac{n^{2+d}}{2} \sum_{w\in \Sn} \int_{Q_n(w)} 
\Bigl(\int_{y\not=x} \bigl(E_nu_n(x)-E_nu_n(y)\bigr)^2 
C_{NN}^n(w,[y]_n)dy \Bigr) dx } \\
& & \vspace{-7pt} \\
&\ge& \dis{ \frac{n^{2+d}}{2} \sum_{w\in \Sn} 
\int_{Q_n(w)}\bigl(E_nu_n(x)-E_nu_n(x+\mbox{\boldmath $e$}_i/n)\bigr)^2 
\Bigl(\int_{Q_n(w+\mbox{\boldmath $e$}_i/n)} dy \Bigr) dx } \\
& & \vspace{-7pt} \\
&=& \dis{ \frac{n^{2}}{2} \sum_{w\in \Sn} 
\int_{Q_n(w)}\bigl(E_nu_n(x)-E_nu_n(x+\mbox{\boldmath $e$}_i/n)\bigr)^2
dx } \\
& & \vspace{-7pt} \\
&=& \dis{ \frac{n^{2}}{2} \int_{\rd}
\bigl(E_nu_n(x)-E_nu_n(x+\mbox{\boldmath $e$}_i/n)\bigr)^2 dx. } \\
\end{array}
$$
In other words, $\{n\bigl(E_nu_n(\cdot)-
E_nu_n(\cdot+\mbox{\boldmath $e$}_i/n)\bigr)\}_n$ is a bounded sequence 
in $L^2(\rd, dx)$.  So there exists a subsequence $\{n'\}$ and a unique 
$v_i\in L^2(\rd, dx)$ so that $n'\bigl(E_{n'}u_{n'}(\cdot)-
E_{n'}u_{n'}(\cdot+\mbox{\boldmath $e$}_i/{n'})\bigr)$ converges to $v_i$ 
{ weakly in} $L^2(\rd, dx)$. On the other hand, if $\varphi 
\in C_c^{2}(\rd)$, it follows that 
$$
\langle E_{n'}u_{n'}(\cdot+\mbox{\boldmath $e$}_i/{n'}), \varphi \rangle
=\langle E_{n'}u_{n'}, \varphi(\cdot-\mbox{\boldmath $e$}_i/{n'}) 
\rangle
$$
by a change of variables, and then 
$$
n'\langle E_{n'}u_{n'}(\cdot+\mbox{\boldmath $e$}_i/n'),\varphi \rangle
-n' \langle E_{n'}u_{n'},\varphi \rangle=
n'\langle E_{n'}u_{n'},\varphi(\cdot-\mbox{\boldmath $e$}_i/{n'}) -\varphi 
\rangle.
$$
Since $\varphi\in C_c^{2}(\rd)$, we see that 
$n'\bigl(\varphi(\cdot-\mbox{\boldmath $e$}_i/{n'}) -\varphi\bigr)$ converges 
to $-\partial\varphi/\partial x_i$ uniformly and in $L^2(\rd, dx)$. 
So we have,  letting  $n'\to\infty$, 
$$
\langle v_i,\varphi \rangle =-\langle H, \partial \varphi/\partial x_i \rangle,
$$
since $u_n$ converges to $H$ uniformly on compact sets.  
This shows that $v_i=\partial H/\partial x_i$ and so $H\in W^{1,2}(\rd)$.

\bigskip
{\it Step 2.} \ We next show that for some subsequence $\{n'\}$, 
$$
\form^{n'}_C(u_{n'},g) \longrightarrow \frac 12 \int_{\rd} \nabla H(x) \cdot 
a(x) \nabla g(x)dx=\form_C(H,g)
$$
for any $g\in C_c^2(\rd)$. Recall (\ref{(6)}); 
since $C_C^n(x,y)=0$ if $|x-y|>\vareps_n$ and the $w,z$ are on the shortest 
paths from $x$ and $y$, it is enough to consider $w$'s only for $|w-z|\le 
\vareps_n$ in the sum of the right hand side of (\ref{(6)}). So
$$
\begin{array}{lcl} \form^n_C(u_n,g)
&=& \dis{ \frac 1{2n^{d}}  \sum_{i,j=1}^d \sum_{z\in \Sn} \nabla_{1/n}^i u_n(z) 
\sum_{w\in \Sn \atop |w-z|\le\vareps_n} \nabla_{1/n}^j g(w) 
{G}_{ij}^n(w,z) } \\

& & \vspace{-6pt} \\
&=& \dis{ \frac1{2n^{d}}  \sum_{i,j=1}^d \sum_{z\in \Sn} \nabla_{1/n}^i u_n(z) 
\nabla_{1/n}^j g(z) \sum_{w\in \Sn \atop |w-z|\le\vareps_n} 
{G}_{ij}^n(w,z) } \\

& & \vspace{-6pt} \\
& & \quad \quad + \ \dis{ \frac 1{2n^{d}}  \sum_{i,j=1}^d \sum_{z\in \Sn} 
\nabla_{1/n}^i u_n(z) \sum_{w\in \Sn \atop |w-z|\le\vareps_n} 
\Bigl(\nabla_{1/n}^j g(w) - \nabla_{1/n}^j g(z) \Bigr) 
{G}_{ij}^n(w,z) } \\

& & \vspace{-6pt} \\
&=:& I^n_1+I^n_2.
\end{array}
$$
Let $K$ be the support of $g\in C_c^2(\rd)$. Since 
$1/n \le \vareps_n\le 1$ and $|w-z|\le \vareps_n$ in the summation defining
$I^n_2$, the  
$z$'s must lie in the set $K_1\cap \Sn$, where $K_1=\{x\in \rd : \ 
 d(K, x)\le 1\}$.
By using the mean value theorem for $g$ and the definition of 
$\nabla_{1/n}^i u_n$, we see that for some 
$0<\theta, \tilde{\theta}<1$ depending on $z$ and $w$, 
$$
\begin{array}{lcl}
2|I^n_2|  &=& 
\dis{ \Bigl| n^{-d}  \sum_{i,j=1}^d \sum_{z\in \Sn} 
\nabla_{1/n}^i u_n(z) \sum_{w\in K_1\cap \Sn \atop |w-z|\le\vareps_n} 
\Bigl(\nabla_{1/n}^j g(w) - \nabla_{1/n}^j g(z) \Bigr) 
{G}_{ij}^n(w,z) \Bigr| } \\

& & \vspace{-6pt} \\
&= & \dis{ \Bigl| n^{1-d}  \sum_{i,j=1}^d \sum_{z\in \Sn} 
\Bigl(u_n(z+\mbox{\boldmath $e$}_i/n) -u_n(z) \Bigr)  } \\

& & \vspace{-6pt} \\
& & \quad \dis{\times \sum_{w\in K_1\cap \Sn \atop |w-z|\le\vareps_n} 
\Bigl( \partial_j g(w+\theta \mbox{\boldmath $e$}_j/n) -
\partial_j g(w+\tilde{\theta}\mbox{\boldmath $e$}_j/n) 
 \Bigr) {G}_{ij}^n(w,z) \Bigr| } \\

& & \vspace{-6pt} \\
&\le& \dis{ \Bigl(\sup_{|z-z'|\le 1/n} \bigl|u_n(z)-u_n(z')\bigr|\Bigr) 
\cdot \sup_{j} ||\partial_{jj} g||_{\infty} \times
\Bigl( n^{-d}  \sum_{i,j=1}^d \sum_{z\in \Sn} 
\sum_{w\in K_1\cap \Sn \atop |w-z|\le\vareps_n}  \bigl|
{G}_{ij}^n(w,z) \bigr| \Bigr) } \\

& & \vspace{-6pt} \\
&=:& \dis{ \Bigl(\sup_{|z-z'|\le 1/n} \bigl|u_n(z)-u_n(z')\bigr|\Bigr) 
\cdot \sup_{j} ||\partial_{jj} g||_{\infty} \times I^n_3. } \\
\end{array}
$$
We now estimate $I^n_3$. Let $K_2=\{x\in \rd : \ d(K_1, x)\le 1\}$. Then, 
$$
\begin{array}{lcl}
I^n_3 &=& \dis{
n^{-d}  \sum_{i,j=1}^d \sum_{z\in \Sn} 
\sum_{w\in K_1\cap \Sn \atop |w-z|\le\vareps_n}  \Bigl|

\sum_{x,y\in \Sn \atop |x-y|\le \vareps_n} 
\Bigl(P^{x,y}(z+\mbox{\boldmath $e$}_i/n, z)-
P^{x,y}(z, z+\mbox{\boldmath $e$}_i/n)\Bigr) } \\

& & \vspace{-6pt} \\
& & \quad \quad \quad \times \dis{ 
\Bigl(P^{x,y}(w+\mbox{\boldmath $e$}_j/n, w)-
P^{x,y}(w, w+\mbox{\boldmath $e$}_j/n)\Bigr) C_C^n(x,y)\Bigr|} \\

& & \vspace{-6pt} \\
&\le& \dis{n^{-d} \sum_{x,y\in \Sn \atop |x-y|\le \vareps_n} 
C_C^n(x,y)  \sum_{i=1}^d \sum_{z\in \Sn}  
\Bigl(P^{x,y}(z+\mbox{\boldmath $e$}_i/n, z) + 
P^{x,y}(z, z+\mbox{\boldmath $e$}_i/n)\Bigr) } \\

& & \vspace{-6pt} \\
& & \quad \quad \quad \times \dis{ 
 \sum_{j=1}^d \sum_{w\in K_1\cap \Sn \atop |w-z|\le\vareps_n} 

\Bigl(P^{x,y}(w+\mbox{\boldmath $e$}_j/n, w)+
P^{x,y}(w, w+\mbox{\boldmath $e$}_j/n)\Bigr). } \\

& & \vspace{-6pt} \\
&=& \dis{n^{-d} \sum_{x \in K_2\cap \Sn, y\in \Sn \atop |x-y|\le \vareps_n} 
C_C^n(x,y)  \sum_{i=1}^d \sum_{z\in \Sn}  
\Bigl(P^{x,y}(z+\mbox{\boldmath $e$}_i/n, z) + 
P^{x,y}(z, z+\mbox{\boldmath $e$}_i/n)\Bigr) } \\

& & \vspace{-6pt} \\
& & \quad \quad \quad \times \dis{ 
 \sum_{j=1}^d \sum_{w\in K_1\cap \Sn \atop |w-z|\le\vareps_n} 
\Bigl(P^{x,y}(w+\mbox{\boldmath $e$}_j/n, w)+
P^{x,y}(w, w+\mbox{\boldmath $e$}_j/n)\Bigr). } \\
\end{array}
$$
The last equality holds since the $w$'s (belonging to $K_1$) lie on 
some shortest path between $x$ and $y$ in the summations for some $x,y\in \Sn$ 
with $|x-y|\le \vareps_n$. Noting now that 
$$
\begin{array}{l}
\hspace{-1cm} \dis{ \sum_{j=1}^d 
\sum_{w\in K_1\cap \Sn \atop |w-z|\le\vareps_n} 
\Bigl(P^{x,y}(w+\mbox{\boldmath $e$}_j/n, w)+ 
P^{x,y}(w, w+\mbox{\boldmath $e$}_j/n)\Bigr) } \\

\vspace{-6pt} \\
\le \dis{ \sum_{j=1}^d 
\sum_{w \in \Sn} 
\Bigl(P^{x,y}(w+\mbox{\boldmath $e$}_j/n, w)+ 
P^{x,y}(w, w+\mbox{\boldmath $e$}_j/n)\Bigr)=n|x-y|} \\
\end{array}
$$
and similarly 
$$
\sum_{i=1}^d 
\sum_{z \in \Sn} 
\Bigl(P^{x,y}(z+\mbox{\boldmath $e$}_i/n, z)+ 
P^{x,y}(z, z+\mbox{\boldmath $e$}_i/n)\Bigr)=n|x-y|,
$$
we see that, using (\ref{cond-bound}),  
$$
I^n_3 \le n^{-d} \sum_{x\in K_2 \cap \Sn} 
\sum_{y\in \Sn \atop |x-y|\le\vareps_n}
n^2|x-y|^2 C_C^n(x,y) \le M \mu^n(K_2),
$$
where $M$ is the constant in the assumption (A3) (see (2.1)). 
So, $I^n_3$ is uniformly bounded in $n$ and hence $I^n_2$ converges to $0$ as $n$ 
tends to $\infty$ since the $\{u_n\}$ are equicontinuous.

Finally we consider the term $I^n_1$:
$$
\begin{array}{lcl}
I^n_1 &=& \dis{ \frac 1{2n^{d}}  \sum_{i,j=1}^d \sum_{z\in \Sn} \nabla_{1/n}^i u_n(z) 
\nabla_{1/n}^j g(z) F^n_{ij}(z) } \\
& & \vspace{-6pt} \\
&=& \dis{ \frac 12\sum_{i,j=1}^d \int_{\rd} \nabla_{1/n}^i E_n u_n(x) \nabla_{1/n}^j 
E_n g(x) F^n_{ij}(x)dx.} \\
\end{array}
$$
Observe that if $f_n$ converges to $f$ weakly in $L^2$ and $g_n$ converges
to $g$ boundedly and almost everywhere, then $f_ng_n$ converges to $fg$ weakly.
To see this, if $h\in L^2$, 
$$\int (f_ng_n)h-\int (fg)h=\int f_n(g_n-g)h+ \Big[\int f_ngh-\int fgh\Big].$$
The    term inside the brackets  on the right hand side  goes to 0 since $f_n$
converges to $f$ weakly and  the boundedness of $g$ implies that $gh$
is in $L^2$. The first term on the right hand side is bounded, using Cauchy-Schwarz, by
$\| f_n\|_2 \,\| (g_n-g)h\|_2$. The factor $\|f_n\|_2$
is uniformly bounded since $f_n$ converges weakly in $L^2$, while
$\|(g_n-g)h\|_2$ converges to 0 by dominated convergence.

Since some subsequence of $\nabla_{1/n}^i E_n u_n$ converges to 
$v_i=\partial_i H$ weakly in $L^2$ (as proved in Step 1), 
and for some further subsequence 
$F_{ij}^n$ converges to $a_{ij}$ boundedly and almost everywhere  
(by (A4) and Remark \ref{RR5.3}) and $\nabla_{1/n}^j g$ converges to $\partial_j g$ 
uniformly on compact sets (because $g\in C_c^2(\R^d)$), 
we see that, along this further   subsequence, the right hand side goes to 
$$
\frac 12 \sum_{i,j=1}^d \int_{\rd} \partial_j H\, \partial_j g\,  a_{ij}\, dx=
\frac 12 \int_{\rd} \nabla H(x) \cdot a(x) \nabla g(x)dx.
$$
Hence 
$$
\form_C^{n'}(u_{n'}, g) \rightarrow \form_C(H,g).
$$
This completes the proof of \eqref{WCElqpw} and hence the theorem. \qed

\begin{flushright}
$\begin{array}{l}
\mbox{Richard F. Bass}\\
\mbox{Department of Mathematics}\\
\mbox{University of Connecticut}\\
\mbox{Storrs, CT 06269, U.S.A.}\\
\mbox{E-mail: {\tt bass@math.uconn.edu}}\\
\mbox{\ }\\
\mbox{Takashi Kumagai}\\
\mbox{Department of Mathematics}\\
\mbox{Kyoto University}\\
\mbox{Kyoto 606-8502, Japan}\\
\mbox{E-mail: {\tt kumagai@math.kyoto-u.ac.jp}}\\
\mbox{\ }\\
\mbox{Toshihiro Uemura}\\
\mbox{School of Business Administration}\\
\mbox{University of Hyogo}\\
\mbox{Kobe, 651-2197, Japan}\\
\mbox{E-mail: {\tt uemura@biz.u-hyogo.ac.jp}}\\
\end{array}$
\end{flushright}
\end{document}